\documentclass[12pt]{amsart}
\usepackage{epsfig,amsmath,amsfonts,latexsym,accents}
\usepackage[usenames,dvipsnames]{color}
\usepackage{palatino}
\usepackage[hypcap=false]{caption}
\usepackage{caption, subcaption, pinlabel}
\usepackage{color}

\usepackage{xcolor}

\definecolor{darkgreen}{rgb}{0, 0.5, 0}

\usepackage[hyperfootnotes=false]{hyperref}
\hypersetup{
  colorlinks,
  citecolor=Purple,
  linkcolor=Black,
  urlcolor=Green}
  
\theoremstyle{plain}
\newtheorem{dummy}{anything}[section]
\newtheorem{theorem}[dummy]{Theorem}

\newtheorem{lemma}[dummy]{Lemma}

\newtheorem{proposition}[dummy]{Proposition}

\theoremstyle{definition}
\newtheorem{definition}[dummy]{Definition}

\newtheorem{remark}[dummy]{Remark}

\theoremstyle{remark}

\newcommand{\Z}{\mathbb{Z}}

\newcommand{\bP}{\mathbb P}

\def\a{\alpha}
\def\b{\beta}

\def\g{\gamma}

\begin{document}

\title {Minimal number of singular fibers in a nonorientable Lefschetz fibration}

\author{Sinem Onaran} 

\author{Burak Ozbagci}

\address{Department of Mathematics, Hacettepe University, 06800, Ankara, Turkey}

\email{sonaran@hacettepe.edu.tr}

\address{Department of Mathematics, Ko\c{c} University, 34450, Istanbul,Turkey}
\email{bozbagci@ku.edu.tr}

\subjclass[2000]{}
\thanks{}


\begin{abstract}

We show that there exists an admissible nonorientable genus $g$ Lefschetz fibration with only one singular fiber over a closed orientable surface of genus $h$ if and only if $g \geq 4$ and $h \geq 1$.

\end{abstract}

\maketitle

\section{Introduction} 

Let $S$ be a nonorientable closed surface and let $t_\a$ denote a Dehn twist about a two-sided simple closed curve $\a$ in $S$. We say that $\a$  is {\em trivial} if it bounds a disk or a M\"{o}bius band.  Note that $\a $ is nontrivial  if and only if $t_\a$ is nontrivial in the mapping class group $\mathcal{M} (S)$  of $S$, which is the group of isotopy classes of all self-diffeomorphisms of $S$.

\begin{definition} Let $X$ be a closed nonorientable $4$-manifold, and $\Sigma$ be a closed orientable $2$-manifold. We say that a (nonorientable) Lefschetz fibration $f: X \to \Sigma$  is {\em admissible} if  it has  at least one singular fiber and all the vanishing cycles of $f$ are nontrivial.  \end{definition}

A Lefschetz fibration is said to be of genus $g$ if the regular fiber is of genus $g$. In a nonorientable Lefschetz fibration,  the regular  fiber is a nonorientable surface, whose genus is defined to be the (maximal) number of projective planes in a connected sum decomposition. 
Note that the interactions between the relations in the mapping class groups of surfaces and Lefschetz fibrations listed in \cite[Section 15.2]{os} for the orientable case hold verbatim for the nonorientable case. In particular, for any closed nonorientable genus $g$ surface $S$, if a product of $m$ nontrivial Dehn twists is equal to a product of $h$ commutators in $\mathcal{M} (S)$,  then there is an admissible (nonorientable) genus $g$ Lefschetz fibration with $m$ singular fibers over a closed orientable surface of genus $h$. 

For any $g \geq 1$ and $h \geq 0$, let $N(g, h)$ denote the minimal number of singular fibers in an admissible  genus $g$ Lefschetz fibration on a closed nonorientable $4$-manifold over a closed orientable surface of genus $h$.  According to \cite[Proposition 1.11]{mo}, any admissible genus one Lefschetz fibration  over any closed orientable surface is a  $\bP^2$-bundle. 
Therefore, $N(1,h)$ is not defined for any $h \geq 0$, since  there are no singular fibers in this case. Moreover, $N(g, 0) >1 $ for all $g \geq 2$, since a Dehn twist along a nontrivial curve is not isotopic to the identity,  by definition. 

 On the other hand, any even power of any Dehn twist on a nonorientable surface $S$ can be expressed as a single commutator in $\mathcal {M} (S)$ (see, for example,  the Remark at the end of  \cite{st}). In particular, for any $g \geq 2$, there exists an admissible  genus $g$ Lefschetz fibration with {\em two} singular fibers over any closed orientable surface of genus $h \geq 0$, which implies that $N(g,h) \leq 2$ if   $g \geq 2$ and $h \geq 0$. As a consequence of this discussion, a natural question arises: For which values of the pair $(g,h)$, we have $N(g, h)=1$?

Our main result is the following 

\begin{theorem} \label{thm: numb} $N(g,h) = 1$ iff $g \geq 4$ and $h \geq 1$. 
\end{theorem}

By our definitions in this paper, Theorem~\ref{thm: numb} is about nonorientable Lefschetz fibrations. An analogous result for the minimal number of singular fibers in {\em orientable} relatively minimal Lefschetz fibrations was obtained by Korkmaz and the second author \cite{ko}, which says that  $N(g,h) = 1$ iff $g \geq 3$ and $h \geq 2$.

Let $[G,G]$ denote the commutator subgroup of a group  $G$. For any $x \in [G,G]$, the commutator length of $x$ is defined to be the minimum number of factors needed to express $x$ as a product of commutators. 
Our proof of Theorem~\ref{thm: numb} is based on the commutator lengths of Dehn twists on nonorientable surfaces.  In \cite{sz}, Szepietowski  showed that the commutator length of any power of any Dehn twist on a closed nonorientable surface of genus  at least $7$ is equal to $1$.  Here, we prove the following results in the remaining low-genus cases. 

\begin{proposition}\label{prop: main} Let $S$ be a closed nonorientable surface of genus  $4, 5$ or $6$. Then for any nontrivial separating curve $\g$ in $S$ and for every $n \in \Z$, the commutator length of $t^n_\g$ is equal to $1$ in  $\mathcal{M} (S)$.
\end{proposition}

\begin{proposition} \label{prop: g3} If $S$ is a closed nonorientable surface of genus $2$ or $3$, then no Dehn twist along a nontrivial curve in $S$ belongs to $[\mathcal{M} (S), \mathcal{M} (S)]$. 
\end{proposition} 

In the rest of the paper, we assume that the genus of a closed nonorientable surface is at least $2$, since  the mapping class group of the real projective plane is trivial. We also assume that all curves are simple closed and two-sided. 

\section{Commutator lengths of Dehn twists on  nonorientable surfaces} 

Let $S$ be a closed nonorientable surface. When $S$ is of genus at least $7$, Szepietowski \cite{sz} proved that the commutator length of every power of any Dehn twist  is  equal to $1$  in $\mathcal{M} (S)$.  He also showed that if $S$ is of genus $6$, then the commutator length of every power of the Dehn twist along any nonseparating curve with orientable complement  is equal to $1$  in $\mathcal{T} (S)$, the twist subgroup of $\mathcal{M} (S)$ which is generated by all Dehn twists. Note that, if the genus of $S$ is at least $7$, then we have $$ [\mathcal{M} (S), \mathcal{M} (S)]= [\mathcal{T} (S), \mathcal{T} (S)]= \mathcal{T} (S)$$  as shown in \cite[Theorem 5.12]{k3}, which is not true  if the genus of $S$ is at most $6$,   since in these cases $H_1( \mathcal{T} (S))$ is nontrivial \cite{st}. We hope that this sheds some light  on the genus assumptions in the aforementioned results of Szepietowski. In this note, we partially extend his results to cover the cases of genus $4,5$ and $6$, which is stated as Proposition~\ref{prop: main} in the Introduction.

\begin{proof} of {\em Proposition~\ref{prop: main}.} Let $S$ be a closed nonorientable surface of genus $4$ and let $\g$ be a nontrivial separating curve in $S$. Note  that at least one component of the complement of $\g$ is  a one-holed Klein bottle.  In that component, we have   $t_\g = y^2$,  where $y$ is the crosscap
transposition. But we also know that $y$ is conjugate to $y^{-1}$ by a diffeomorphism $f$ of the surface $S$. It follows that  $$ t^n_\g=y^ny^n=y^nfy^{-n}f^{-1}. $$
 
Next,  let $S$ be a closed nonorientable surface of genus $5$ or $6$ and let $\g$ be a nontrivial separating curve in $S$. When we cut the surface $S$ along $\g$, we see that at least one of the components must be {\em nonorientable} of genus $2, 3$ or $4$. If there is a nonorientable component of genus  $2$ in the complement of $\g$, then we can prove the desired result exactly as described in the first paragraph. Now suppose that one of the nonorientable components in the complement of $\g$ is of genus  $3$. Recall that there is the following well-known relation 
\begin{equation} \label{eq:20}
t_d t_e = (t_a t_b t_c)^4
\end{equation} on the two-holed torus,  where $d$ and $e$ are the boundary components.  There is an orientation-reversing involution $r$ of the two-holed torus such that $r (a) =c$,   $r(b) = b$, and $r(c) =a$. Since r is orientation
reversing,  it conjugates right-handed Dehn twists to left-handed Dehn twists, and as a result we have 
\begin{equation} \label{eq:21}
r(t_a t_b t_c)^{-2}r =(t_a t_b t_c)^2
\end{equation}
using that fact $r= r^{-1}$.  By inserting the relation in (\ref{eq:21}) into  (\ref{eq:20}), we obtain
\begin{equation} \label{eq:31} 
t_d t_e = r(t_a t_bt_c)^{-2}r(t_a t_b t_c)^2.  \end{equation} 
Finally we embed the two-holed torus into our surface $S$ such that the boundary component $d$ is capped off with a  M\"{o}bius band, and $e$ is identified with $\g$.
Then $t_d$ becomes trivial in (\ref{eq:31}) and hence $t_\g$ is a single commutator in $\mathcal{M} (S)$.  Note that the involution $r$ extends onto $S$. It follows  that, for all $n \in \Z$,  the commutator length of $t_\g^n$ is also equal to $1$  as well. 

The only remaining case is that $\g$ separates a closed nonorientable surface $S$ of genus $6$ into a one-holed torus and a nonorientable surface of genus $4$. To take care of this case, we first make some preliminary observations. In a one-holed torus, let $c$ denote a boundary parallel curve and let $a, b$ denote the standard generators of the first homology group. The relation $(t_a t_b)^6 =t_c $ is well-known in the mapping class group of the one-holed torus.  Using the braid relation $t_a t_b t_a = t_b t_a t_b$ we get 
\begin{equation} \label{eq:1}
 t_c = (t_bt_a^2t_bt_a^2)^2. 
\end{equation}
 There is an orientation-reversing involution $r$ of the one-holed torus such that $r (a) =a$ and $r(b) = b$. Since r is orientation
reversing,  it conjugates right-handed Dehn twists to left-handed Dehn twists, and as a result we have 
\begin{equation} \label{eq:2}
r(t_bt_a^2t_bt_a^2)^{-1}r = t_a^2t_bt_a^2t_b.
\end{equation}
Here we used that fact $r= r^{-1}$.  By using the braid relation again we get 
\begin{equation} \label{eq:3} 
t_a^{-1} (t_a^2t_bt_a^2t_b) t_a  = t_bt_a^2t_bt_a^2.
\end{equation}
Combining (\ref{eq:2}) and (\ref{eq:3}) we obtain 
\begin{equation} \label{eq:4} 
t_a^{-1}r (t_bt_a^2t_bt_a^2)^{-1} rt_a = t_bt_a^2t_bt_a^2.
\end{equation}
Finally, we insert the relation in (\ref{eq:4}) into the relation in (\ref{eq:1}), to get  
\begin{equation} \label{eq:5} 
t_c = (t_bt_a^2t_bt_a^2)^2 = t_a^{-1}r (t_bt_a^2t_bt_a^2)^{-1} rt_a (t_bt_a^2t_bt_a^2) = [t_bt_a^2t_bt_a^2,   t_a^{-1}r]. \end{equation} Similarly, for any $n \in \Z$, we have 
\begin{equation} \label{eq:6} 
t^n_c = (t_bt_a^2t_bt_a^2)^n (t_bt_a^2t_bt_a^2)^n = t_a^{-1}r (t_bt_a^2t_bt_a^2)^{-n} rt_a (t_bt_a^2t_bt_a^2)^{n} = [(t_bt_a^2t_bt_a^2)^n,   t_a^{-1}r]. \end{equation} There is an obvious embedding of the one holed torus  into $S$  so that the boundary $c$ is identified with $\g$  in $S$ and the involution $r$ of the one holed torus extends to an involution of $S$. This shows, combined with (\ref{eq:6}),  that any power of the Dehn twist $t_\g$ can be expressed as a single commutator in $\mathcal{M} (S)$. 
\end{proof}

Next we show that one cannot possibly extend Szepietowski's results to cover the cases of genus $2$ or $3$, which is stated as Proposition~\ref{prop: g3} in the Introduction.

\begin{proof} of {\em Proposition~\ref{prop: g3}.} The mapping class group of the closed nonorientable surface $S$ of genus $2$  is isomorphic to $\Z_2 \oplus \Z_2$ (cf. \cite{l})  and therefore its commutator subgroup is trivial.  Hence the Dehn twist along the only nontrivial curve in $S$ cannot be equal to a product of commutators. 

Suppose that  $S$ is a closed nonorientable surface of genus $3$ for the rest of the proof. It follows from the presentation of  $\mathcal{M} (S)$ given in \cite{bc} that the homology group $H_1 (\mathcal{M} (S)) = \Z_2 \oplus \Z_2$ is generated by the homology class $\bar{t}_{a}$ of a Dehn twist  $t_{a}$ and the homology class $\bar{z}$  of a crosscap slide $z$. In fact (cf. \cite{k}) we have $$H_1 (\mathcal{M} (S)) = \langle \bar{t}_{a}, \bar{z} :  {\bar{t}}^2_{a}= \bar{z}^2=1, \bar{t}_{a} \bar{z} = \bar{z} \bar{t}_{a}  \rangle .$$   Suppose that $b$ is a nonseparating curve in $S$ with nonorientable complement. Then $t_{b}$ is conjugate to $t_{a}$. Since conjugate elements in a group $G$ are homologous in $H_1(G)$, we have $\bar{t}_{b} = \bar{t}_{a}$. We conclude that  $t_{b} \notin [\mathcal{M} (S), \mathcal{M} (S)]$,  since otherwise, $\bar{t}_{a} = \bar{t}_{b}=1$, which indeed contradicts to the fact that  $\bar{t}_{a}$ is a generator of $H_1 (\mathcal{M} (S))$.  Note that there is no nonseparating curve with {\em orientable} complement in $S$.  Moreover, since any separating curve in $S$ is trivial, we deduce that no Dehn twist along a nontrivial curve in $S$ is a commutator.
\end{proof} 
 
 Using similar techniques, one can prove the following result.

\begin{lemma} \label{lem: g26} If $S$ is a closed nonorientable surface of genus $4, 5$ or $6$, then no Dehn twist along a nonseparating curve with nonorientable complement in $S$ belongs to $[\mathcal{M} (S), \mathcal{M} (S)]$. 
\end{lemma} 

\begin{remark} Notice that there is no such curve if $S$ is of genus $2$ and the genus $3$ case is already covered by Proposition~\ref{prop: g3}. Moreover, if $S$ is of genus $4$, then it follows from \cite{k} that the Dehn twist about
a nonseparting curve with {\em orientable} complement represents a nontrivial element of $H_1 (\mathcal{M} (S))$, and hence it does not belong to $[\mathcal{M} (S), \mathcal{M} (S)]$.  \end{remark} 

We now turn to the proof of our main Theorem~\ref{thm: numb}.  

\begin{proof} {\em Theorem~\ref{thm: numb}.} Szepietowski \cite[Theorem 1.1]{sz} showed that if $S$ is a closed nonorientable surface of genus at least  $7$, then for every nontrivial curve $\a$  in $S$, $t_\a$  is equal to a single commutator in $\mathcal{M} (S)$. Moreover, he showed (cf. \cite[Theorem 1.2]{sz}) that there exists a nontrivial curve $\b$ in  a closed nonorientable surface of genus $6$, such that  $t_{\b}$ is equal to a single commutator of elements of  the twist subgroup $\mathcal{T} (S) \subset \mathcal{M} (S)$. It follows that for any $g \geq 6$ and $h \geq 1$, there is  an admissible  genus $g$ Lefschetz fibration over a closed orientable surface of genus $h$, which has a unique singular fiber. The cases $g=4,5$ follow from Proposition~\ref{prop: main}. The cases $g=2,3$ follow from Proposition~\ref{prop: g3}. \end{proof}

\section{A Lemma in the orientable case} 
In this final section, $S$ denotes a closed  {\em orientable} surface of positive genus.  As shown by Korkmaz and the second author \cite{ko}, if $S$ is of genus at least $3$, then the commutator length of any  Dehn twist   is equal to $2$ in $\mathcal{M} (S)$. This result is optimum and there is no room for improvement. This is because,  for a closed orientable surface $S$ of genus $1$ or $2$, no Dehn twist belongs to $[\mathcal{M} (S), \mathcal{M} (S)]$ (see, for example, \cite{ko}) since $H_1 (\mathcal{M} (S))$ is nontrivial in these cases. It is also known that, for  any closed orientable surface $S$,  the commutator length of any power of a Dehn twist about a separating curve is at least $2$ (see \cite[Remark 3]{ek}) and the same is true for a sufficiently high power of a Dehn twist about a nonseparating curve \cite[Corollary 2.4]{k2}.    

The extended mapping class group of an orientable surface $S$, denoted by $\mathcal{M}^\diamond (S)$, includes also the isotopy classes of orientation-reversing self-diffeomorphisms of $S$. In \cite{sz}, Szepietowski showed that if $S$ is of genus at least $3$, then the commutator length of every power of any Dehn twist  is equal to $1$ in $\mathcal{M}^\diamond (S)$.  It turns out that there is some room for improvement for his result. 

\begin{lemma}\label{lem: two}  Let $\g$ be a nontrivial separating curve in a closed orientable surface $S$ of genus $2$. Then for every $n \in \Z$, the commutator length of $t^n_\g$ is equal to $1$ in  $\mathcal{M}^\diamond (S)$. 
\end{lemma} 

\begin{remark} Note that for a closed orientable surface $S$ of genus $1$ or $2$, no Dehn twist along a nonseparating curve belongs to $[\mathcal{M}^\diamond (S), \mathcal{M}^\diamond  (S)]$. This follows from the fact that all nonseparating curves on a fixed surface $S$ are equivalent and the homology class of the Dehn twist along some nonseparating curve is nontrivial in $H_1 ( \mathcal{M}^\diamond (S)) = \Z_2 \oplus \Z_2$ (\cite[Theorem 4.2]{k}).  We conclude that Lemma~\ref{lem: two} cannot be improved any further. \end{remark}

\begin{proof}{\em of Lemma~\ref{lem: two}.}    Let $\g$ be a nontrivial separating curve in a closed orientable surface $S$ of genus $2$. There is an obvious embedding of the one holed torus  into $S$ so that the boundary $c$ is identified with $\g$  in $S$ and the involution $r$ of the one holed torus extends to an involution of $S$. This shows, combined with (\ref{eq:6}),  that any power of the Dehn twist $t_\g$ can be expressed as a single commutator in $\mathcal{M}^\diamond (S)$. \end{proof}

\noindent {\bf {Acknowledgement}}: We are grateful to Mustafa Korkmaz and the  referee whose suggestions helped us improve the results in this paper. SO was partially supported by the Young Scientist Awards Program BAGEP of the Science Academy, Turkey.  \\


\end{document}